\documentclass[11pt]{amsart}
\usepackage[margin=2.5cm]{geometry}
\usepackage{enumerate}
\usepackage{amsmath}
\usepackage{amssymb,latexsym}
\usepackage{amsthm}
\usepackage{color}
\usepackage[dvipsnames]{xcolor}
\usepackage{cancel}
\usepackage{listings}
\usepackage{graphicx}
\usepackage{cite}
\usepackage{changes}
\usepackage{mathdots}
\usepackage{float}
\usepackage{hyperref}
\usepackage[normalem]{ulem}

\makeatletter
\renewcommand{\tocsection}[3]{%
  \indentlabel{\@ifnotempty{#2}{\bfseries\ignorespaces#1 #2\quad}}\bfseries#3}
\renewcommand{\tocsubsection}[3]{%
  \indentlabel{\@ifnotempty{#2}{\ignorespaces#1 #2\quad}}#3}

\renewcommand{\l@paragraph}[4]{%
  \@tocline{4}{0pt}{1.5em}%
  {\indentlabel{\@ifnotempty{#1}{\ignorespaces#1\quad}}}%
  {}%
}
\makeatother

\usepackage{etoolbox}
\makeatletter
\patchcmd{\@setaddresses}{\indent}{\noindent}{}{}
\patchcmd{\@setaddresses}{\indent}{\noindent}{}{}
\patchcmd{\@setaddresses}{\indent}{\noindent}{}{}
\patchcmd{\@setaddresses}{\indent}{\noindent}{}{}
\makeatother

\newcommand{\F}{\mathbb F}

\newcommand{\rk}{\operatorname{rk}}

\newcommand{\Tr}{\operatorname{Tr}}

\newcommand{\PG}{\mathrm{PG}}

\newcommand{\Sym}{\mathrm{Sym}}

\newcommand{\rad}{\operatorname{rad}}

\usepackage{cleveref}

\theoremstyle{definition}
\newtheorem{theorem}{Theorem}[section]
\newtheorem{problem}[theorem]{Problem}
\newtheorem{lemma}[theorem]{Lemma}
\newtheorem{corollary}[theorem]{Corollary}
\newtheorem{definition}[theorem]{Definition}
\newtheorem{proposition}[theorem]{Proposition}

\newtheorem{remark}[theorem]{Remark}

\title{Non-additive symmetric rank-metric codes}
\date{}

\author[M. Messia]{Marta Messia}
\address{Marta Messia, \textnormal{Department of Mathematics and Data Science, Vrije Universiteit Brussel, Pleinlaan 2, 1050 Brussel, Belgium}}
\email{marta.messia@vub.be}

\author[F. Zullo]{Ferdinando Zullo}
\address{Ferdinando Zullo, \textnormal{Dipartimento di Matematica e Fisica, Universit\`a degli Studi della Campania ``Luigi Vanvitelli'', Viale Lincoln, 5, I--\,81100 Caserta, Italy}}
\email{ferdinando.zullo@unicampania.it}

\begin{document}

\begin{abstract}
We study non-additive rank-metric codes contained in the space of symmetric matrices over a finite field. Our main focus is on symmetric $2$-codes. For odd $q$, we introduce a quotient construction for $\mathrm{Sym}_3(q)$ in which rank-one differences are described by isotropic vectors of a nondegenerate quadratic form on $\mathbb F_q^4$. This leads to a correspondence with partial ovoids of the elliptic quadric $Q^-(5,q)$. Using a construction of Cossidente and Pavese, we obtain a non-additive symmetric $2$-code in $\mathrm{Sym}_3(q)$ of cardinality $q^4+2q^3$, improving the previously known general lower bound. We also give an explicit matrix description of this family. Lastly, we present a block-based construction for symmetric $2$-codes that extends the three-dimensional improvement to all odd dimensions. In particular, for every odd $n\geq3$ and every prime power $q>2$, we obtain
\begin{equation*}
    A_q^{\rm sym}(n,2) \ge q^{(n^2-1)/2}+q^{(n^2-3)/2}+q^{(n^2-9)/2},
\end{equation*}
which is strictly larger than the maximum cardinality of an additive symmetric $2$-code with the same parameters. This yields the first infinite family (as the matrix dimension $n$ varies) of non-additive symmetric $2$-codes whose size is larger than the largest achievable size of an additive code with the same parameters; before this, such a gap had only been established for $n=3$.
\end{abstract}

\maketitle

\tableofcontents

\section{Introduction}
\label{sec:introduction}

Rank-metric codes were introduced by Delsarte in the framework of bilinear forms over finite fields~\cite{Delsarte1978}.  Given two matrices $A$ and $B$ over $\mathbb F_q$, their rank distance $d_r(A,B)$ is $\operatorname{rk}(A-B)$. Besides the classical setting of arbitrary rectangular matrices, it is natural to consider rank-metric codes contained in matrix spaces with additional algebraic structure, such as symmetric, alternating, or Hermitian matrices. These spaces arise naturally in finite geometry and association schemes, and the corresponding extremal problems may behave quite differently from those in the unrestricted matrix space; \cite{Delsarte1978,Schmidt2015,SchmidtHermitian2018,CossidenteMarinoPavese2022}.

In this paper, we are concerned with symmetric matrices.  Let $\mathrm{Sym}_n(q)=\{A\in M_n(\mathbb F_q)\,:\,A=A^T\}$. A subset
\[
\mathcal C\subseteq \mathrm{Sym}_n(q)
\]
is called a \emph{symmetric $d$-code} if
\[
\operatorname{rk}(A-B)\ge d\qquad\text{for all distinct }A,B\in\mathcal C.
\]
We denote by $A_q^{\mathrm{sym}}(n,d)$ the largest possible cardinality of a symmetric $d$-code in $\mathrm{Sym}_n(q)$.

A particularly important distinction is whether or not the code is assumed to be additive. If $\mathcal C$ is an additive subgroup of $\mathrm{Sym}_n(q)$, then we call it an \emph{additive symmetric rank-metric code}. We write $M_{\mathrm{add}}^{\mathrm{sym}}(n,d)$ for the largest cardinality of an additive symmetric $d$-code.

The additive theory of symmetric rank-metric codes is well understood through the association schemes of symmetric bilinear forms; see Schmidt~\cite{Schmidt2015,Schmidt2020}.  In particular, when $d$ is even, then we have
\begin{equation}
A_q^{\mathrm{sym}}(n,d)\le\begin{cases}
q^{\,n(n-d+2)/2},
& n-d \text{ even},\\[2mm]
q^{\,(n+1)(n-d+1)/2},
& n-d \text{ odd}.
\end{cases}
\label{eq:symmetric-additive-bound-intro}
\end{equation}
When $d$ is odd, the same bounds hold for a non additive code. These bounds are attained by suitable additive codes. In particular, for odd $n$,
\begin{equation}
M_{\mathrm{add}}^{\mathrm{sym}}(n,2)=q^{(n^2-1)/2}.
\label{eq:additive-symmetric-2}
\end{equation}
The construction and classification of maximum additive symmetric rank-metric codes remain an interesting area of research; see, for instance, \cite{TangZhou2026,MushrrafZullo2025,AlnajjarineLavrauw2026,LongobardiLunardonTrombettiZhou2020,Zhou2020}.

The situation changes substantially when additivity is removed. In the case of $d$ even, the bounds in \eqref{eq:symmetric-additive-bound-intro} are bounds for additive codes and need not hold for arbitrary subsets of $\mathrm{Sym}_n(q)$. This leads to the basic extremal problem of how large can a non-additive symmetric rank-metric code be. In particular, one may ask when the additive optimum can be exceeded, and by how much.

Cossidente, Marino and Pavese~\cite{CossidenteMarinoPavese2022} proved that, for $q>2$, there exists a symmetric $2$-code in $\mathrm{Sym}_3(q)$ of cardinality $q^4+q^3+1$, which is strictly larger than the additive optimum $q^4$; cf. \eqref{eq:additive-symmetric-2}. They also obtained the upper bound
\begin{equation}
A_q^{\mathrm{sym}}(3,2)\le 1+\frac{q(q^2-1)(q^2+q+1)}{2}.
\label{eq:CMP-sym-upper}
\end{equation}

The first aim of this paper is to improve the known general lower bound for this problem.  For every odd prime power $q$, we construct a symmetric $2$-code in $\mathrm{Sym}_3(q)$ of cardinality $q^4+2q^3$. Consequently,
\[
A_q^{\mathrm{sym}}(3,2)\ge q^4+2q^3,
\]
and the codes obtained in this way are necessarily non-additive.

Our construction is based on a geometric interpretation of a suitable quotient of the space of symmetric matrices. More precisely, for odd $q$ we consider a two-dimensional subspace of $\mathrm{Sym}_3(q)$ whose nonzero elements all have rank two, and study the quotient of $\mathrm{Sym}_3(q)$ by this subspace. The key observation is that, in the quotient, rank-one differences are described by the zero set of a nondegenerate quadratic form on a four-dimensional vector space.  In this way, the
problem of constructing a symmetric $2$-code is reduced to finding a large set whose nonzero pairwise differences avoid this quadric.

The latter problem has a natural interpretation in finite polar geometry. After embedding the four-dimensional quadratic space into a six-dimensional space of minus type, the relevant quadratic form gives rise to the elliptic quadric $Q^-(5,q)$.  Under a suitable affine parameterization, the condition that two vectors have an anisotropic difference corresponds precisely to the associated points of $Q^-(5,q)$ being nonorthogonal, or equivalently noncollinear. Consequently, partial ovoids of $Q^-(5,q)$ provide exactly the geometric configurations needed for our construction.

This establishes a direct connection between partial ovoids and non-additive symmetric rank-metric codes. In particular, a partial ovoid of $Q^-(5,q)$ of size $m$ gives rise to a symmetric $2$-code in $\mathrm{Sym}_3(q)$ of size $q^2(m-1)$. Thus, results on large partial ovoids can be translated directly into lower bounds for $A_q^{\mathrm{sym}}(3,2)$.

Cossidente and Pavese~\cite{CossidentePavese2014} constructed, for every prime power $q$, a maximal partial ovoid of $Q^-(5,q)$ of size $(q+1)^2$.  For odd $q$, their construction, therefore, yields

$$
A_q^{\mathrm{sym}}(3,2)\ge q^4+2q^3,
$$

improving the previously known general construction.  We also make this correspondence explicit by parameterizing the Cossidente--Pavese partial ovoid into a concrete family of symmetric matrices.

More generally, the extensive literature on partial ovoids of $Q^-(5,q)$, see for instance \cite{DeBeuleKleinMetschStorme2008,Coolsaet2014}, provides further input for our approach: any improvement in the size of such geometric configurations immediately produces a corresponding improvement for symmetric rank-metric codes.

A second feature of our construction is that an improvement in dimension three can be propagated to arbitrarily large odd dimensions. For symmetric $2$-codes $\mathcal C\subseteq\mathrm{Sym}_{n_1}(q)$ and $\mathcal D\subseteq\mathrm{Sym}_{n_2}(q)$, we introduce the block product
\[
\mathcal C\star\mathcal D=\left\{\begin{pmatrix}
A&X\\
X^T&B
\end{pmatrix}\,:\,A\in\mathcal C,\;B\in\mathcal D,\;X\in\mathbb F_q^{\,{n_1}\times {n_2}}\right\}.
\]
We show that $\mathcal C\star\mathcal D$ is again a symmetric $2$-code and that
\begin{equation*}
|\mathcal C\star\mathcal D|=q^{n_1n_2}|\mathcal C||\mathcal D|,
\end{equation*}
combining this construction with a two-dimensional additive construction propagates any three-dimensional symmetric $2$-code to every odd dimension.

In particular, starting from the partial-ovoid construction above, we obtain, for every odd prime power $q$ and every odd integer $n\ge3$,
\begin{equation*}
A_q^{\mathrm{sym}}(n,2)\ge q^{(n^2-1)/2} + 2q^{(n^2-3)/2}.
\end{equation*}
Moreover, for every prime power $q>2$ and every odd integer $n\ge3$ we show 
    \begin{equation*}
        A_q^{\rm sym}(n,2) \ge q^{(n^2-1)/2}+q^{(n^2-3)/2}+q^{(n^2-9)/2}.
    \end{equation*}
Comparing with \eqref{eq:additive-symmetric-2}, we exhibit two infinite families of symmetric $2$-codes which are strictly larger than the largest possible additive codes with the same parameters.

The paper is organized as follows. In Section~2, we recall the basic notions and results on symmetric rank-metric codes, quadratic spaces, and partial ovoids that will be used throughout the paper. In Section~3, we introduce the quotient construction for $\Sym_3(q)$ and show that rank-one differences in the quotient are described by isotropic vectors of a nondegenerate quadratic form on $\F_q^4$. In Section~4, we interpret this quadratic condition geometrically in terms of partial ovoids of the elliptic quadric $Q^-(5,q)$, obtaining our main lower bound for $A_q^{\rm sym}(3,2)$. Section~5 is devoted to the Cossidente--Pavese construction: we determine it explicitly in the affine coordinates arising from our quotient model and derive a concrete family of non-additive symmetric $2$-codes of cardinality $q^4+2q^3$. In Section~6, we introduce a block construction that allows symmetric $2$-codes to be combined and use it to propagate the three-dimensional construction to all odd dimensions. Finally, Section~7 contains some concluding remarks and discusses several open problems suggested by our results.

\section{Preliminaries}
\label{sec:preliminaries}

In this section, we collect the notation and the basic facts on symmetric rank-metric codes and quadrics over finite fields that will be used throughout the paper. We keep the background to a minimum and refer to \cite{Schmidt2015,Schmidt2020} for symmetric rank-metric codes and to \cite{HirschfeldThas} for the standard theory of quadrics over finite fields.

\subsection{Symmetric rank-metric codes}

Let $q$ be a prime power and let
$\mathrm{Sym}_n(q)=\{A\in M_n(\mathbb F_q)\,:\, A=A^T\}$ be the space of $n\times n$ symmetric matrices over $\mathbb F_q$. The \emph{rank distance} on $\mathrm{Sym}_n(q)$ is defined by $d_r(A,B)=\operatorname{rk}(A-B)$. A subset $\mathcal C\subseteq\mathrm{Sym}_n(q)$ is called a \emph{symmetric $d$-code} if
\[
d_r(A,B)\ge d\qquad\text{for all distinct } A,B\in\mathcal C.
\]
We denote by $A_q^{\mathrm{sym}}(n,d)$ the maximum cardinality of a symmetric $d$-code in $\mathrm{Sym}_n(q)$.

A symmetric $d$-code $\mathcal C$ is called \emph{additive} if it is an additive subgroup of $\mathrm{Sym}_n(q)$.  We denote by $M_{\mathrm{add}}^{\mathrm{sym}}(n,d)$ the maximum cardinality of an additive symmetric $d$-code.

The following result summarizes the bounds that will be relevant for us. It follows from the association-scheme approach to symmetric bilinear forms developed by Schmidt; see \cite{Schmidt2015} for odd $q$ and \cite[Theorem~4.3]{Schmidt2020} for the general statement.

\begin{theorem}
\label{prop:Schmidt-bound}
Let $\mathcal C\subseteq\mathrm{Sym}_n(q)$ be a symmetric $d$-code. If $d$ is odd, then
\begin{equation*}
|\mathcal C|\le\begin{cases}
q^{\,n(n-d+2)/2},& n-d \text{ even},\\[1mm]
q^{\,(n+1)(n-d+1)/2},& n-d \text{ odd}.
\end{cases}
\end{equation*}
If $d$ is even, the same bounds hold under the additional assumption that $\mathcal C$ is additive.
\end{theorem}

For even minimum distance, the bound is, in general, only an additive bound, and non-additive symmetric codes may be larger. For odd $n$,
\begin{equation*}
M_{\mathrm{add}}^{\mathrm{sym}}(n,2)=q^{(n^2-1)/2}.
\end{equation*}
As discussed in the introduction, the fact that the additive bound can be exceeded when $d$ is even is precisely the phenomenon that motivates our constructions.

\subsection{Quadratic spaces and partial ovoids}

We now recall the geometric terminology needed in the sequel.  Since our constructions will be over fields of odd order, throughout this subsection we assume that $q$ is odd.

Let $V$ be a finite-dimensional vector space over $\mathbb F_q$, and let $Q:V\longrightarrow\mathbb F_q$ be a quadratic form.  Its associated symmetric bilinear form is
\[
\beta_Q(x,y)=\frac{Q(x+y)-Q(x)-Q(y)}{2}.
\]
Thus
\[
Q(x+y)=Q(x)+Q(y)+2\beta_Q(x,y).
\]
Two vectors $x,y\in V$ are said to be \emph{orthogonal}, and we write $x\perp y$, if $\beta_Q(x,y)=0$. The quadratic form $Q$ is \emph{nondegenerate} if 
\begin{equation*}
    \rad(\beta_Q)=\{x\in V\,:\,\beta_Q(x,y)=0\text{ for every }y\in V\}=\{0\}.
\end{equation*}

For a nonzero vector $x\in V$, we write $\langle x\rangle$ for the corresponding point of the projective space $\mathrm{PG}(V)=\mathrm{PG}(\dim(V)-1,q)$. If $V$ has dimension six and $Q$ is a nondegenerate quadratic form of minus type, the set of points defined by vectors that annihilate $Q$ form an \emph{elliptic quadric}, denoted by $Q^-(5,q)$. We refer to \cite[Chapter~1]{HirschfeldThas} for the classification and basic properties of nondegenerate quadrics over finite fields.

We shall repeatedly use the standard relation between orthogonality and collinearity on a quadric. If $P=\langle x\rangle,R=\langle y\rangle$ are two distinct points of a nondegenerate quadric, then 
\begin{equation}
\langle P,R\rangle\subseteq Q\quad\Longleftrightarrow\quad\beta_Q(x,y)=0.
\label{eq:collinear-orthogonal}
\end{equation}
Thus the whole projective line joining $P$ and $R$ is contained in the quadric precisely when $P$ and $R$ are orthogonal. In this case, we will write $P \perp R$.

Let $Q:V\to\mathbb F_q$ be a nondegenerate quadratic form.  The \emph{projective general orthogonal group} associated with $Q$ is the subgroup of $\mathrm{PGL}(V)$ induced by the linear transformations $g\in\mathrm{GL}(V)$ for which there exists $\lambda_g\in\mathbb F_q^*$ such that

$$
Q(gx)=\lambda_g Q(x)\qquad\text{for every }x\in V.
$$

We denote this group by $\mathrm{PGO}(Q)$.  Since such transformations preserve the zero set of $Q$, the group $\mathrm{PGO}(Q)$ acts naturally on the corresponding projective quadric.

We recall also a standard consequence of Witt's theorem.
\begin{theorem}
\label{th:acttrans}
Let $Q$ be a nondegenerate quadratic form, then $\mathrm{PGO}(Q)$ acts transitively on the points of the corresponding projective quadric.
\end{theorem}

Thus, given any two points $P$ and $R$ of the quadric, there exists an element of $\mathrm{PGO}(Q)$ mapping $P$ to $R$; see \cite{HirschfeldThas}. We refer to \cite{HirschfeldThas} for the standard properties of the projective orthogonal groups and their action on finite quadrics.

\begin{definition}
A \emph{partial ovoid} of $Q^-(5,q)$ is a set $\mathcal O$ of points of the quadric such that every line contained in $Q^-(5,q)$ meets $\mathcal O$ in at most one point.
\end{definition}
\begin{remark}
Observe that, by \eqref{eq:collinear-orthogonal}, a set $\mathcal O\subseteq Q^-(5,q)$ is a partial ovoid if and only if no two distinct points of $\mathcal O$ are orthogonal.   
\end{remark}

A partial ovoid is called \emph{maximal} if it is not properly contained in another partial ovoid.  We reserve the term \emph{maximum} for a partial ovoid of largest possible cardinality, and write
\[
\mu^-(5,q)=\max\bigl\{|\mathcal O|\,:\,\mathcal O \text{ is a partial ovoid of }Q^-(5,q)\bigr\}.
\]

\begin{remark}
\label{rem:preserve-ovoid}
The definition of a partial ovoid relies on the relation of collinearity between the points of the considered quadric. Since every element of $\mathrm{PGO}(Q)$ preserves this relation, partial ovoids are preserved under the action of $\mathrm{PGO}(Q)$.
\end{remark}

Specific constructions and bounds for partial ovoids will be recalled later when they are applied to the construction of symmetric rank-metric codes.

\section{The quotient construction}

Throughout this section, assume that $q$ is odd.  Fix a nonsquare $\delta\in\F_q^*$ and define
\begin{equation*}
H_\delta=\left\{\begin{pmatrix}
a&b&0\\
b&\delta a&0\\
0&0&0
\end{pmatrix}\,:\,a,b\in\F_q\right\}\leq\Sym_3(q).
\end{equation*}

\begin{lemma}
\label{lem:Hdelta-rank}
Let $q$ be odd. Every nonzero element of $H_\delta$ has rank $2$.
\end{lemma}

\begin{proof}
The upper-left $2\times2$ determinant is $\delta a^2-b^2$. If it vanished for $(a,b)\ne(0,0)$, then either $a=0$ and hence $b=0$, or $(b/a)^2=\delta$, contradicting that $\delta$ is nonsquare. Thus, the $2\times2$ block is nonsingular whenever $(a,b)\ne(0,0)$, while the final row and column are zero.
\end{proof}

Since $H_\delta$ is a subgroup, we will consider elements in $\Sym_3(q)/H_\delta$: each coset is a symmetric $2$-code, because the difference of two distinct elements belongs to $H_\delta$, and therefore has rank $2$. 

Now, we will identify $H_\delta$ with the kernel of a homomorphism, that will help us to represent all the matrices whose difference has rank $1$. Define the linear map
\begin{align*}
\pi:\qquad\Sym_3(q)\quad&\rightarrow\quad\quad\quad\F_q^4\\
\begin{pmatrix}
x&y&z\\
y&w&t\\
z&t&u
\end{pmatrix}&\mapsto(w-\delta x,z,t,u).
\end{align*}

Clearly, $\pi$ is surjective and $\ker\pi=H_\delta$, hence, we get 
\[
\Sym_3(q)/H_\delta\simeq\F_q^4.
\]

On $\F_q^4$, define the following quadratic form
\begin{equation}
Q(R,Z,T,U)=T^2-\delta Z^2-RU.
\label{eq:Qminus}
\end{equation}

The next lemma is the key point of the construction.

\begin{lemma}
\label{lem:rank-one-quotient}
Let $q$ be odd. If  $A\in\Sym_3(q)$ has rank one then $Q(\pi(A))=0$. Conversely, if $v \in \F_q^4\setminus\{0\}$ is such that $Q(v)=0$, then there exists a rank-one matrix $A\in\Sym_3(q)$ such that $\pi(A)=v$. Equivalently,
\[
\pi\bigl(\{A\in\Sym_3(q)\,:\,\rk A=1\}\bigr)=\{v\in\F_q^4\setminus\{0\}\,:\,Q(v)=0\}.
\]
\end{lemma}

\begin{proof}
Let $A\in\Sym_3(q)$ with $\rk A=1$, then there exist $c\in\F_q^*$ and some nonzero vector $w=(p,r,s)^T\in\F_q^3$ such that
\[
A=c\,ww^T=c\begin{pmatrix}
p^2&pr&ps\\
pr&r^2&rs\\
ps&rs&s^2
\end{pmatrix}.
\]
Hence, $\pi(A)=\bigl(c(r^2-\delta p^2),\,cps,\,crs,\,cs^2\bigr)$, substituting into \eqref{eq:Qminus} gives
\[
Q(\pi(A))=c^2r^2s^2-\delta c^2p^2s^2-c(r^2-\delta p^2)\,cs^2=0.
\]
Moreover, since every matrix in $\ker\pi$ has rank $2$, $\pi(A)\neq0$ for every rank-$1$ matrix $A$. Therefore, $\pi\bigl(\{A\in\Sym_3(q):\rk A=1\}\bigr)\subseteq\{v\in\F_q^4\setminus\{0\}:Q(v)=0\}$.

Conversely, let $(R,Z,T,U)\in\F_q^4$ be a nonzero vector satisfying $Q((R,Z,T,U))=T^2-\delta Z^2-RU=0$. 

Suppose first that $U\neq0$ and set
\[
c=U,\qquad w=\left(\frac ZU,\frac TU,1\right).
\]
Then we can define the rank-one symmetric matrix
\[
A:=c\,ww^T=\begin{pmatrix}
\frac {Z^2}{U}&\frac {ZT}{U}&Z\\
\frac {ZT}{U}&\frac {T^2}{U}&T\\
Z&T&U
\end{pmatrix}.
\]
If we apply the map $\pi$ we obtain
\begin{equation*}
    \pi(A)=\left(\frac{T^2-\delta Z^2}{U},Z,T,U\right)=(R,Z,T,U),
\end{equation*}
since the vector $(R,Z,T,U)$ satisfies the quadratic form $Q$. Therefore, $(R,Z,T,U)$ is the image under $\pi$ of the rank-one matrix $A$.

Suppose now that $U=0$, from the quadratic equation we get $T^2=\delta Z^2$, as $\delta$ is a nonsquare, this forces $T=Z=0$. Given that $(R,Z,T,U)$ is nonzero, $R\neq0$ and we can set
\[
c=R,\qquad w=\left(0,1,0\right).
\]
Following the same reasoning,
\[
A:=c\,ww^T=\begin{pmatrix}
0&0&0\\
0&R&0\\
0&0&0
\end{pmatrix},
\]
applying the map $\pi$ we get $\pi(A)=(R,0,0,0)$. Therefore, $\pi\bigl(\{A\in\Sym_3(q):\rk A=1\}\bigr)\supseteq\{v\in\F_q^4\setminus\{0\}:Q(v)=0\}$ and the result is proved.
\end{proof}

The following result allows us to restrict the problem of building a symmetric $2$-code in $\Sym_3(q)$ to that of finding a set $S\subseteq\F_q^4$ such that the difference of any two distinct elements is not a zero of the quadratic form.

\begin{proposition}
\label{prop:quotient-construction}
Let $q$ be odd and let $S\subseteq\F_q^4$ satisfy
\begin{equation}
Q(x-y)\neq0\qquad\text{for all distinct }x,y\in S.
\label{eq:anisotropic-differences}
\end{equation}
Then $\mathcal{C}_S:=\pi^{-1}(S)$ is a symmetric $2$-code in $\Sym_3(q)$, and
\begin{equation*}
|\mathcal{C}_S|=q^2|S|.
\end{equation*}
\end{proposition}

\begin{proof}
Since $|\ker\pi|=|H_\delta|=q^2$, the preimage of every element of $S$ under $\pi$ contains exactly $q^2$ matrices. Therefore $|\mathcal C_S|=q^2|S|$. Now we want to prove that $\mathcal C_S$ has minimum distance $2$. Let $A,B\in\mathcal C_S$ be distinct and suppose that $\pi(A)=\pi(B)$, then $A-B\in H_\delta\setminus\{0\}$, and Lemma~\ref{lem:Hdelta-rank} gives $\rk(A-B)=2$.

Suppose instead that $\pi(A)\neq\pi(B)$, and, by linearity, $\pi(A-B)=\pi(A)-\pi(B)$. If $\rk(A-B)=1$, Lemma~\ref{lem:rank-one-quotient} would imply $Q\bigl(\pi(A)-\pi(B)\bigr)=0$, contrary to \eqref{eq:anisotropic-differences}.  Hence $\rk(A-B)\ge2$, i.e., $\mathcal C_S$ is a $2$-code.
\end{proof}

When the code is the union of cosets of $H_{\delta}$, the above result is a characterization.

\begin{proposition}
\label{prop:quotient-construction iff}
Let $q$ be odd. Let $\mathcal{S}$ be a union of cosets of $H_\delta$, $\mathcal{S}$ is a symmetric $2$-code if and only if 
\begin{equation*}
    Q(x-y)\neq0\qquad\text{for all distinct }x,y\in \pi (\mathcal{S}).
\end{equation*}
\end{proposition}
\begin{proof}
Clearly, by Proposition \ref{prop:quotient-construction} if $Q(x-y)\neq0$ for all distinct $x,y\in \pi (\mathcal{S})$, then $\mathcal{S}$ is a symmetric $2$-code.

Conversely, suppose that $\mathcal{S}$ is a symmetric $2$-code and, by contradiction, that there exist distinct $x,y\in\pi(\mathcal{S})$ such that $Q(x-y)=0$. Since $x,y\in\pi(\mathcal{S})$, there exist $A,B\in\mathcal{S}$ such that $\pi(A)=x$ and $\pi(B)=y$. Moreover, given $Q(x-y)=0$, by Lemma~\ref{lem:rank-one-quotient} and the linearity of $\pi$, there exists $C\in\Sym_3(q)$ such that $\rk(C)=1$ and
\begin{equation*}
    \pi(C)=x-y=\pi(A-B)=\pi(A)-\pi(B).
\end{equation*} 

Suppose first that $A$ and $B$ belong to the same coset, thus $A-B\in H_\delta$. Since $\pi(C)=\pi(A-B)$ we have that $A-B-C\in H_\delta$, hence $C\in H_\delta$. However, $\rk(C)=1$ and every matrix in $H_\delta$ has rank equal to $2$, therefore a contradiction.

Now, take $A$ and $B$ in two different cosets and let $D:=C+B$, then $\pi(D)=\pi(C)+\pi(B)=\pi(A)$, so $D$ and $A$ belong to the same coset, thus $D\in\mathcal{S}$, since $\mathcal{S}$ is a $2$-code then $2\leq\rk(D-B)=\rk(C)=1$, again a contradiction.
\end{proof}

\section{Partial ovoids and symmetric rank-metric codes}

In this section, we will construct a large set $S\subseteq\F_q^4$ with the property that
\begin{equation*}
Q(x-y)\neq0\qquad\text{for all distinct }x,y\in S,
\end{equation*}
where
\[
Q(R,Z,T,U)=T^2-\delta Z^2-RU.
\]
Recall that such a set gives, by Proposition~\ref{prop:quotient-construction}, a symmetric $2$-code in $\Sym_3(q)$ of size $q^2|S|$.

Let $x\in\F_q^4$ and $a,b\in\F_q$, then we will write the elements of $\F_q^6$ as $(a,x,b)$ and define
\begin{equation*}
\widetilde Q(a,x,b)=ab-Q(x).
\end{equation*}

\begin{lemma}
\label{lem:Qtilde-minus-type}
Let $q$ be odd and let $\delta\in\F_q^*$ be a nonsquare.  The quadratic form
\[
\widetilde Q(a,R,Z,T,U,b)=ab+RU+\delta Z^2-T^2
\]
is nondegenerate and its projective points form an elliptic quadric $Q^-(5,q)$.
\end{lemma}

\begin{proof}

Firstly, we observe that $\widetilde{Q}$ is nondegenerate. Indeed, with respect to the ordered coordinates $(a,R,Z,T,U,b)$, the symmetric matrix associated with $\widetilde{Q}$ is

$$
M_{\widetilde{Q}}
=
\frac{1}{2}
\begin{pmatrix}
0 & 0 & 0 & 0 & 0 & 1\\
0 & 0 & 0 & 0 & 1 & 0\\
0 & 0 & 2\delta & 0 & 0 & 0\\
0 & 0 & 0 & -2 & 0 & 0\\
0 & 1 & 0 & 0 & 0 & 0\\
1 & 0 & 0 & 0 & 0 & 0
\end{pmatrix}.
$$

Since $q$ is odd and $\delta\neq0$, we have $\det(M_{\widetilde{Q}})
=
-\frac{\delta}{16}\neq0$.

Thus $M_{\widetilde{Q}}$ is nonsingular, and consequently $\widetilde{Q}$ is a nondegenerate quadratic form.

Set $f(Z,T)=\delta Z^2-T^2$, since $\delta$ is a nonsquare in $\mathbb{F}_q$, the binary quadratic form $f$ is irreducible over $\mathbb{F}_q$.

Finally, we can write
\begin{equation*}
   \widetilde{Q}(a,R,Z,T,U,b)=f(Z,T)+RU+ab, 
\end{equation*}

where $f(Z,T)$ is an irreducible binary quadratic form over $\mathbb{F}_q$. Hence $\widetilde{Q}$ has the canonical form of a non-degenerate elliptic quadric, see \cite[Section 1.1]{HirschfeldThas}. Therefore, its projective points form an elliptic quadric $Q^-(5,q)$.
\end{proof}

Consider now the vector 
\begin{equation*}
    e:=(1,0,0,0,0,0).
\end{equation*}

Since $\widetilde Q(e)=0$, the associated projective point $P_\infty=\langle e\rangle$ belongs to the quadric $Q^-(5,q)$.

The next lemma gives an affine parameterization for the points that are not orthogonal to $P_\infty$.

\begin{lemma}
\label{lem:affine-chart}
Let $q$ be odd. The points of $Q^-(5,q)$ that are not orthogonal to $P_\infty=\langle (1,0,0,0,0,0)\rangle$ are in one-to-one correspondence with the vectors of $\F_q^4$. More precisely, the correspondence is
\begin{equation*}
x\in \F_q^4\longmapsto P_x=\left\langle\bigl(Q(x),x,1\bigr)\right\rangle \in Q^-(5,q).
\end{equation*}
Moreover, for $x,y\in\F_q^4$,
\begin{equation}
P_x\perp P_y\quad\Longleftrightarrow\quad Q(x-y)=0.
\label{eq:orthogonality-difference}
\end{equation}
\end{lemma}

\begin{proof}
Write a vector of $\F_q^6$ in the form
\[
X=(a,x,b),\qquad x=(R,Z,T,U)\in\F_q^4.
\]
Then $\widetilde Q(X)=ab-Q(x)$.

We first determine which points of the quadric are not orthogonal to $P_\infty$. The polar form of $\widetilde Q$ is
\[
\beta_{\widetilde Q}\bigl((a,x,b),(a',y,b')\bigr)=\frac{ab'+a'b}{2}-\beta_Q(x,y).
\]
In particular, let $e=(1,0,0,0,0,0)$,
\[
\beta_{\widetilde Q}(e,(a,x,b))=\frac b2\neq0\quad\Longleftrightarrow\quad b\neq0.
\]

Therefore, a projective point $\langle(a,x,b)\rangle$ is not orthogonal to $P_\infty$ precisely when $b\neq0$. Now, let $\langle(a,x,b)\rangle$ be a point of the quadric not orthogonal to $P_\infty$. Since $b\neq0$, we can choose a representative of the form $(c,x,1)$. Given that $\langle(c,x,1)\rangle$ belongs to the quadric, we have $0=\widetilde Q(c,x,1)=c-Q(x)$. Therefore, the point is of the form
\[
P_x=\langle(Q(x),x,1)\rangle,
\]
this proves the claimed bijection with $\F_q^4$.

It remains to determine when two such points are orthogonal. Using the polar form above,
\[
\beta_{\widetilde Q}\bigl((Q(x),x,1),(Q(y),y,1)\bigr)=\frac{Q(x)+Q(y)}{2}-\beta_Q(x,y).
\]
On the other hand,
$Q(x-y)=Q(x)+Q(y)-2\beta_Q(x,y)$.
Hence
\[
\beta_{\widetilde Q}((Q(x),x,1),(Q(y),y,1))=\frac12Q(x-y).
\]
Since $q$ is odd, $\beta_{\widetilde Q}((Q(x),x,1),(Q(y),y,1))=0$ if and only if $Q(x-y)=0$ and we conclude the proof.
\end{proof}

Thus, our original problem has the following geometric reformulation: we seek a large collection of points $P_x$ that belongs to the quadric associated to $\widetilde Q$ no two of which are orthogonal, i.e. a partial ovoid in $Q^-(5,q)$.

\begin{theorem}
\label{prop:partial-ovoid-input}
Let $q$ be odd. Suppose that $Q^-(5,q)$ contains a partial ovoid $\mathcal O$ of size $m$. Then there exists a set $S\subseteq\F_q^4$ of size $|S|=m-1$, such that
\[
Q(x-y)\neq0\qquad\text{for every distinct }x,y\in S.
\]
Consequently,
\begin{equation*}
A_q^{\rm sym}(3,2)\ge q^2(m-1).
\end{equation*}
\end{theorem}

\begin{proof}
Let $P_0\in\mathcal O$, by Theorem \ref{th:acttrans} the orthogonal group of $Q^-(5,q)$ acts transitively on its points. Since every element of this group maps a partial ovoid to a partial ovoid (see Remark \ref{rem:preserve-ovoid}), we may therefore apply such a transformation and assume, without loss of generality, that $P_0=P_\infty=\langle(1,0,0,0,0,0)\rangle$.

Since $\mathcal O$ is a partial ovoid, no point of $\mathcal O\setminus\{P_\infty\}$ is orthogonal to $P_\infty$. Hence, Lemma~\ref{lem:affine-chart} applies to every such point. Consequently, for every $\bar{P}\in\mathcal O\setminus\{P_\infty\}$ there is a unique vector $x\in\F_q^4$, such that
\[
\bar{P}=P_x=\langle (Q(x),x,1)\rangle.
\]

Define $S=\left\{x\in\F_q^4\,:\,P_x\in\mathcal O\setminus\{P_\infty\}\right\}$. Since the correspondence $x\mapsto P_x$ is one-to-one,
\[
|S|=|\mathcal O|-1=m-1.
\]

Now, let $x,y\in S$ be distinct. Then $P_x$ and $P_y$ are two distinct points of the partial ovoid $\mathcal O$.  By definition of a partial ovoid, $P_x\not\perp P_y$, therefore, by Lemma~\ref{lem:affine-chart}, we have $Q(x-y)\neq0$, hence, $S$ satisfies \eqref{eq:anisotropic-differences}.

Finally, Proposition~\ref{prop:quotient-construction} shows that $\pi^{-1}(S)\subseteq\Sym_3(q)$ is a symmetric $2$-code of size $q^2|S|=q^2(m-1)$, thus
\[
A_q^{\rm sym}(3,2)\ge q^2(m-1).
\]
\end{proof}

We can now apply a known construction of partial ovoid. Cossidente and Pavese in \cite{CossidentePavese2014} constructed, for every prime power $q$, a partial ovoid of $Q^-(5,q)$ having size $(q+1)^2$. For odd $q$, combining their construction with Theorem~\ref{prop:partial-ovoid-input} gives the following family of symmetric rank-metric codes.

\begin{corollary}
\label{cor:main-symmetric}
For every odd prime power $q$,
\begin{equation*}
A_q^{\rm sym}(3,2)\ge q^4+2q^3.
\end{equation*}
Moreover, the code obtained in this way is non-additive.
\end{corollary}

\begin{proof}
Take a partial ovoid of size $m=(q+1)^2$, whose existence is guaranteed by \cite[Theorem 6.1]{CossidentePavese2014}. Theorem~\ref{prop:partial-ovoid-input} produces a set $S\subseteq\F_q^4$ of size
\[
|S|=m-1=(q+1)^2-1=q^2+2q.
\]
The quotient construction then gives a symmetric $2$-code of size
\[
q^2|S|=q^2(q^2+2q)=q^4+2q^3.
\]

The maximum size of an additive symmetric $2$-code in $\Sym_3(q)$ is $q^4$. Since $q^4+2q^3>q^4$, the code constructed above cannot be additive.
\end{proof}

Recall that $\mu^-(5,q)$ denotes the maximum cardinality of a partial ovoid of $Q^-(5,q)$.

\begin{proposition}
\label{prop:fibre-partial-ovoid-equivalence}
Let $q$ be odd. The maximum size of a symmetric $2$-code in $\Sym_3(q)$ which is a union of complete cosets of $H_\delta$ is
\begin{equation*}
    q^2\bigl(\mu^-(5,q)-1\bigr).
\end{equation*}
In particular,
\[ 
A_q^{\rm sym}(3,2)\ge q^2\bigl(\mu^-(5,q)-1\bigr).
\]
\end{proposition}

\begin{proof}
By Proposition~\ref{prop:quotient-construction iff}, such a code is of the form $\pi^{-1}(S)$, where $S\subseteq\F_q^4$ satisfies $Q(x-y)\neq0$ for all distinct $x,y\in S$ and has size $q^2|S|$. For every such set $S$, consider $\mathcal O_S=\{P_\infty\}\cup\{P_x:x\in S\}$. The point $P_\infty$ is not orthogonal to any $P_x$, while Lemma~\ref{lem:affine-chart} gives $P_x\not\perp P_y$ if and only if $Q(x-y)\neq0$. Hence, $\mathcal O_S$ is a partial ovoid of $Q^-(5,q)$ of size $|S|+1$. It follows that $|S|\le\mu^-(5,q)-1$.

Conversely, let $\mathcal O$ be a partial ovoid of maximum size $\mu^-(5,q)$ and choose $P_0\in\mathcal O$.  By transitivity of the orthogonal group on the points of the quadric, we may assume that $P_0=P_\infty$. Every other point of $\mathcal O$ is uniquely of the form $P_x$, and the corresponding set $S=\{x:P_x\in\mathcal O\setminus\{P_\infty\}\}$ has size $|S|=\mu^-(5,q)-1$ and satisfies $Q(x-y)\neq0$ for distinct $x,y\in S$.

Therefore, the maximum possible size of $S$ is $\mu^-(5,q)-1$, and the statement follows from $|\pi^{-1}(S)|=q^2|S|$.
\end{proof}

\begin{remark}
\label{rem:known-partial-ovoids}
More precise information is available for small values of $q$. In particular, Lansdown~\cite[Section~4.7]{Lansdown2020} collects and improves the known upper and lower bounds for the maximum size $\mu^-(5,q)$ of a partial ovoid of $Q^-(5,q)$. Restricting to the odd values relevant to our construction, Theorem~\ref{prop:partial-ovoid-input} gives the following bounds.
\[
\begin{array}{c|c|c|c}
q & \text{bounds on }\mu^-(5,q)
  & \text{resulting lower bound on }A_q^{\rm sym}(3,2) &\text{CMP upper bound}\\
\hline
3  & 16                       & 135 & 157\\
5  & 48                       & 1175 & 1861\\
7  & 98\leq\mu^-(5,7)\leq176 & 4753 & 9577\\
9  & 148\leq\mu^-(5,9)\leq361& 11907 & 32761\\
11 & 258\leq\mu^-(5,11)\leq672&31097 & 87781\\
13 & 330\leq\mu^-(5,13)\leq1106&55601 & 199837
\end{array}
\]

The last column of the table refers to \cite[Corollary~4.15]{CossidenteMarinoPavese2022}, where the authors proved that the cardinality of a symmetric $2$-code in $\Sym_3(q)$ is upper-bounded by the value
\begin{equation*}
    1+\frac{q(q^2-1)(q^2+q+1)}{2}.
\end{equation*}

For $q=3$ and $q=5$ the values of $\mu^-(5,q)$ are exactly known. In particular, Lansdown gives an exhaustive enumeration for $Q^-(5,5)$ and obtains $\mu^-(5,5)=48$; moreover, the maximum partial ovoid is unique up to equivalence. Consequently, Proposition~\ref{prop:fibre-partial-ovoid-equivalence} shows that the maximum cardinality of a symmetric $2$-code obtained from the union of cosets of $H_\delta$ construction is $135$ for $q=3$ and $1175$ for $q=5$.

For the remaining field sizes, the exact value of $\mu^-(5,q)$ is not known. It is worth noting that Lansdown's computations improve the previously known lower bounds, for $q=9,11,13$ from $146,216,273$ to $148,258,330$, respectively. Thus, improvements in the partial-ovoid problem translate immediately, through our correspondence, into improved lower bounds for symmetric rank-metric codes.
\end{remark}

\section{An explicit family from the Cossidente--Pavese construction}
\label{sec:explicit-CP}

In the previous section, we showed how partial ovoids of $Q^-(5,q)$ produce symmetric $2$-codes.  We now make this correspondence explicit for the partial ovoid constructed by Cossidente and Pavese in \cite[Theorem~6.1]{CossidentePavese2014}.

The purpose of this section is to describe their construction in the affine coordinates introduced in Lemma~\ref{lem:affine-chart}. To do this, we first fix the notation for the current section, let $q>3$ be odd, let $\delta\in\F_q^*$ be a nonsquare, for $\eta\in\F_q\setminus\{0,\pm1\}$ define $\lambda:=1-\eta^2$. We define the set $S$ to be
\begin{equation*}
    S=S_0\cup S_1\cup S_2\cup S_3\subseteq\F_q^4,
\end{equation*}
where
\begin{equation}
\begin{aligned}
S_0&=\left\{\left(r,z,-\frac12,\frac{1-4\delta z^2}{4r}\right)\,:\,r\in\F_q^*,\ z\in\F_q\right\},\\
S_1&=\left\{(0,0,t,0)\,:\,t\in\F_q\setminus\{0,-1\}\right\},\\
S_2&=\left\{\left(-\frac{\lambda}{1-\delta t^2},-\frac{\lambda t}{1-\delta t^2},-1,\frac{\lambda\delta t^2}{1-\delta t^2}\right)\,:\,t\in\F_q\right\}\cup\{(0,0,-1,-\lambda)\},\\
S_3&=\left\{\left(\frac{1}{1-\delta t^2},\frac{t}{1-\delta t^2},0,-\frac{\delta t^2}{1-\delta t^2}\right)\,:\,t\in\F_q\right\}\cup\{(0,0,0,1)\}.
\end{aligned}
\label{eq:setS}
\end{equation}
In the following, we will prove that $S$ is a parameterization of a partial ovoid in the construction of Cossidente--Pavese and give an explicit description of the corresponding symmetric rank-metric code. Therefore, the aim of this section is to prove the following result.

\begin{theorem}
\label{thm:explicit-CP-code}
Let $q>3$ be odd, let $\delta\in\F_q^*$ be a nonsquare, and choose $\eta\in\F_q\setminus\{0,\pm1\}$ and  $\lambda=1-\eta^2$.
Let
\[
S=S_0\cup S_1\cup S_2\cup S_3\subseteq\F_q^4
\]
be the set defined by \eqref{eq:setS}. Then
\begin{equation*}
\mathcal C_{\delta,\lambda}=\left\{\begin{pmatrix}
a & b & Z\\
b & \delta a+R & T\\
Z & T & U
\end{pmatrix}\,:\,a,b\in\F_q,\;(R,Z,T,U)\in S\right\}\subseteq\Sym_3(q)
\end{equation*}
is a symmetric $2$-code with cardinality
\[
|\mathcal C_{\delta,\lambda}|=q^2(q^2+2q)=q^4+2q^3.
\]
In particular,
\[
A_q^{\rm sym}(3,2)\ge q^4+2q^3.
\]
Moreover, $\mathcal C_{\delta,\lambda}$ is non-additive.
\end{theorem}

In order to prove Theorem \ref{thm:explicit-CP-code}, we first recall the Cossidente--Pavese construction. Let $\mathcal Q=Q^-(5,q)$ be an elliptic quadric. The construction starts from a solid $\Sigma\simeq\PG(3,q)$ such that $\Sigma\cap\mathcal Q\simeq Q^-(3,q)$. Its polar line $\ell=\Sigma^\perp$ is a secant line of $\mathcal Q$.  We denote the two points of $\ell\cap\mathcal Q$ by $P_1$ and $P_2$.

Then choose a point $P\in\Sigma\setminus \mathcal Q$ and consider the plane $\pi=P^\perp\cap\Sigma$. The intersection $\mathcal C=\pi\cap\mathcal Q$ is a nonsingular conic contained in $\Sigma\cap\mathcal Q$.  Moreover, the polar plane $\pi^\perp$ meets $\mathcal Q$ in another nonsingular conic $\mathcal C'=\pi^\perp\cap\mathcal Q$. The first part of the partial ovoid is
\begin{equation*}
\mathcal O_0 = \bigl((\Sigma\cap\mathcal Q)\setminus\mathcal C\bigr) \cup \bigl(\mathcal C'\setminus\ell\bigr).
\end{equation*}

To complete the construction, choose a second solid $\Sigma'\subseteq P^\perp$ such that $\Sigma'\cap\mathcal Q\simeq Q^-(3,q)$ and satisfying two externality conditions. In particular, both $\Sigma'\cap\ell$ and $\Sigma'\cap\pi$ should be external to $\mathcal Q$\footnote{In fact, their external condition concerns a ${Q}(4,q)$ contained in $\mathcal{Q}$. Hence, the condition stated in this paper is equivalent to the one used in the proof of \cite[Theorem~6.1]{CossidentePavese2014}.}. For every $c\in\mathcal C$, the two generators $\langle P_1,c\rangle$ and $\langle P_2,c\rangle$ meet $\Sigma'$ in two more points.  The set $\mathcal B$ formed by these $2(q+1)$ points is then added to $\mathcal O_0$, and
\[
\Omega=\mathcal O_0\cup\mathcal B
\]
is a partial ovoid of $\mathcal Q$ of size $(q+1)^2$; see \cite[Theorem~6.1]{CossidentePavese2014}.

We now realize this configuration in the elliptic quadric defined in the previous section. Namely, recall the quadratic form
\[
Q(R,Z,T,U)=T^2-\delta Z^2-RU
\]
and the elliptic quadric
\begin{equation}
\mathcal Q=\left\{\langle (a,R,Z,T,U,b)\rangle\in\PG(5,q)\,:\,ab-Q(R,Z,T,U)=0\right\}.
\label{eq:CP-quadric}
\end{equation}

\begin{proposition}
\label{prop:CP-affine-form}
Let $q>3$ be odd, let $\delta\in\F_q^*$ be a nonsquare, and choose $\eta\in\F_q\setminus\{0,\pm1\}$ and $\lambda=1-\eta^2$. Then the Cossidente--Pavese construction yields a partial ovoid
\[
\Omega=\{P_\infty\}\cup\{P_x\,:\,x\in S\}\subseteq\mathcal Q
\]
of size $(q+1)^2$, where
\[
P_\infty=\langle(1,0,0,0,0,0)\rangle
\]
and
\[
P_x=\langle (Q(x),R,Z,T,U,1)\rangle,\qquad x=(R,Z,T,U)\in S.
\]

Where
\[
S=S_0\cup S_1\cup S_2\cup S_3\subseteq\F_q^4
\]
is the set defined by \eqref{eq:setS}.

In particular,
\[
|S|=q^2+2q
\]
and
\[
Q(x-y)\neq0\qquad \text{for all distinct }x,y\in S.
\]
\end{proposition}

\begin{proof}
We parameterize a partial ovoid in the construction of Cossidente--Pavese~\cite[Theorem~6.1]{CossidentePavese2014} in the coordinate system of \eqref{eq:CP-quadric}.

Consider the solid
\[
\Sigma: \begin{cases} a=0, \\ 
b=-2T. \end{cases}
\]
Its intersection with $\mathcal Q$ is
\[
\Sigma\cap\mathcal Q=\left\{\langle(0,R,Z,T,U,-2T)\rangle\,:\,RU+\delta Z^2-T^2=0\right\},
\]
which is an elliptic quadric $Q^-(3,q)$.

Using the polar form of the quadratic form defining $\mathcal Q$, one obtains
\[
\ell=\Sigma^\perp\,:\,\begin{cases}
R=Z=U=0,\\ 
T=-a.
\end{cases}
\]
Intersecting with $\mathcal Q$ gives
\[
\ell\cap\mathcal Q=\{P_1,P_2\},
\]
where
\[
P_1=\langle(1,0,0,-1,0,1)\rangle,\qquad P_2=\langle(0,0,0,0,0,1)\rangle.
\]

Choose
\[
P=\langle(0,0,0,1,0,-2)\rangle.
\]
Since $P\in\Sigma\setminus\mathcal Q$, it can be used as the external point in the Cossidente--Pavese construction.  A direct computation gives
\[
P^\perp\,:\qquad a+T=0.
\]
Hence,
\[
\pi=P^\perp\cap\Sigma\,:\qquad a=b=T=0,
\]
and therefore
\begin{equation}
\mathcal C=\pi\cap\mathcal Q=\left\{\langle(0,R,Z,0,U,0)\rangle\,:\,RU+\delta Z^2=0\right\}.
\label{eq:CP-conic}
\end{equation}

Moreover,
\[
\pi^\perp\,:\qquad R=Z=U=0,
\]
so
\begin{equation*}
\mathcal C'=\pi^\perp\cap\mathcal Q=\left\{\langle (a,0,0,T,0,b)\rangle\,:\,ab=T^2\right\}.
\end{equation*}
Thus, the set $\mathcal O_0=\bigl((\Sigma\cap\mathcal Q)\setminus\mathcal C\bigr) \cup \bigl(\mathcal C'\setminus\ell\bigr)$ is completely determined in our coordinates.

We first describe its points in the affine coordinates introduced in Lemma~\ref{lem:affine-chart}.

A point of $(\Sigma\cap\mathcal Q)\setminus\mathcal C$ has $T\neq0$, and hence, given $q$ odd, its last coordinate $b=-2T$ is nonzero. After normalizing $b=1$, we obtain $T=-\frac12$. Since $a=0$, the equation of $\mathcal Q$ becomes
\[
RU+\delta Z^2-\frac14=0.
\]
Notice that $R\neq0$, otherwise $1=4\delta Z^2$, which is impossible since $\delta$ is a nonsquare. 
Therefore,
\[
U=\frac{1-4\delta Z^2}{4R}.
\]
It follows that $(\Sigma\cap\mathcal Q)\setminus\mathcal C=\{P_x\,:\,x\in S_0\}$.

We next consider $\mathcal C'\setminus\ell$.  The point $P_\infty$ belongs to this set.  Every other point has $b\neq0$, and after normalizing $b=1$ can be written as
\[
\langle (t^2,0,0,t,0,1)\rangle,\qquad t\in\F_q.
\]
Such a point belongs to $\ell$ precisely when $t=-t^2$, that is, when $t=0$ or $t=-1$. Hence,
\[
\mathcal C'\setminus\ell=\{P_\infty\}\cup\{P_x\,:\,x\in S_1\}.
\]
Consequently,
$\mathcal O_0=\{P_\infty\}\cup\{P_x\,:\,x\in S_0\cup S_1\}$.

It remains to determine the additional $2(q+1)$ points forming $\mathcal B$.

We first parametrize the conic $\mathcal C$.  From \eqref{eq:CP-conic},
\[
\mathcal C=\{c_t\,:\,t\in\F_q\}\cup\{c_\infty\},
\]
where
\[
c_t=\langle(0,1,t,0,-\delta t^2,0)\rangle,\qquad t\in\F_q,\qquad c_\infty=\langle(0,0,0,0,1,0)\rangle.
\]

We now explicitly choose the auxiliary solid $\Sigma'$ appearing in the Cossidente--Pavese construction.  Put
\begin{equation*}
\Sigma'_\lambda\,:\,\begin{cases}
T=-a,\\
b=(1+\lambda)a+R+U.
\end{cases}
\end{equation*}
Since $T=-a$, we have $\Sigma'_\lambda\subseteq P^\perp$. Let us briefly check that $\Sigma'_\lambda$ satisfies the externality conditions required. Restricting the quadratic form of $\mathcal Q$ to $\Sigma'_\lambda$ gives
\[
\begin{aligned}
ab+RU+\delta Z^2-T^2&=a\bigl((1+\lambda)a+R+U\bigr)+RU+\delta Z^2-a^2\\
&=\lambda a^2+a(R+U)+RU+\delta Z^2.
\end{aligned}
\]
Setting $R'=R+a$ and $U'=U+a$, this becomes
\[
R'U'+(\lambda-1)a^2+\delta Z^2.
\]
Since $\lambda-1=-\eta^2$, we obtain
\[
R'U'+\delta Z^2-(\eta a)^2,
\]
which is a quadratic form of minus type; see \cite[Section 1.1]{HirschfeldThas}. Hence,
\[
\Sigma'_\lambda\cap\mathcal Q\simeq Q^-(3,q).
\]

Furthermore, evaluating the quadratic form defining $\mathcal Q$ at the point 
\[
\Sigma'_\lambda\cap\ell=\langle(1,0,0,-1,0,1+\lambda)\rangle,
\]
we obtain the value $\lambda$, which is nonzero. Thus, $\Sigma'_\lambda\cap\ell$ is external to $\mathcal Q$.

Similarly,
\[
\Sigma'_\lambda\cap\pi\,:\,\begin{cases}
a=b=T=0,\\
R+U=0.
\end{cases}
\]
If this line meets $\mathcal C$, then $-R^2+\delta Z^2=0$, which has no projective solution because $\delta$ is a nonsquare. Thus, $\Sigma'_\lambda$ satisfies the required externality conditions.

For $c=\langle(0,R,Z,0,U,0)\rangle\in\mathcal C$, consider the two generators
\[
L_1(c)=\langle P_1,c\rangle,\qquad L_2(c)=\langle P_2,c\rangle.
\]
Now we determine their intersections with $\Sigma'_\lambda$. Put
\[
s=R+U.
\]
We have $s\neq0$. Indeed, otherwise $U=-R$, and the equation of $\mathcal C$ would give $-R^2+\delta Z^2=0$, again impossible since $\delta$ is a nonsquare.

A generic point of $L_1(c)$ has the form
\[
\alpha P_1+\beta c=\langle(\alpha,\beta R,\beta Z,-\alpha,\beta U,\alpha)\rangle.
\]
The equation $T=-a$ of $\Sigma'_\lambda$ is automatically satisfied, while $b=(1+\lambda)a+R+U$ gives
\[
\lambda\alpha+\beta s=0.
\]
Taking $\alpha=1$, we obtain
\begin{equation*}
B_1(c)=L_1(c)\cap\Sigma'_\lambda=\left\langle\left(1,-\frac{\lambda R}{s} ,-\frac{\lambda Z}{s} ,-1,-\frac{\lambda U}{s},1\right)\right\rangle.
\end{equation*}

Similarly, a generic point of $L_2(c)$ has the form
\[
\alpha P_2+\beta c=\langle(0,\beta R,\beta Z,0,\beta U,\alpha)\rangle,
\]
again, the equation $T=-a$ is automatically satisfied. The defining equation of $\Sigma'_\lambda$ gives
\[
\alpha=\beta s.
\]
Taking $\alpha=1$, we obtain
\begin{equation*}
B_2(c)=L_2(c)\cap\Sigma'_\lambda=\left\langle\left(0,\frac{R}{s},\frac{Z}{s},0,\frac{U}{s},1\right)\right\rangle.
\end{equation*}

We now substitute the parameterization of $\mathcal C$.

For a point $c_t=\langle(0,1,t,0,-\delta t^2,0)\rangle$ we have $s=1-\delta t^2$. Since $\delta$ is a nonsquare,
\[
1-\delta t^2\neq0\qquad\text{for every }t\in\F_q.
\]
Substituting in $B_1(c_t)$ gives
\[
\left\langle\left(1,-\frac{\lambda}{1-\delta t^2},-\frac{\lambda t}{1-\delta t^2},-1,\frac{\lambda\delta t^2}{1-\delta t^2},1\right)\right\rangle.
\]
The point corresponding to $c_\infty$ gives
\[
\langle(1,0,0,-1,-\lambda,1)\rangle.
\]
Hence the points $B_1(c)$ correspond precisely to the vectors in $S_2$.

Likewise, substituting in $B_2(c_t)$ gives
\[
\left\langle\left(0,\frac{1}{1-\delta t^2},\frac{t}{1-\delta t^2},0,-\frac{\delta t^2}{1-\delta t^2},1\right)\right\rangle,
\]
while $c_\infty$ gives
\[
\langle(0,0,0,0,1,1)\rangle.
\]
Thus the points $B_2(c)$ correspond precisely to the vectors in $S_3$.

We have therefore obtained
\[
\Omega =\{P_\infty\}\cup\{P_x\,:\,x\in S_0\cup S_1\cup S_2\cup S_3\},
\]
and this is exactly the Cossidente--Pavese partial ovoid.

Finally,
\[
|S_0|=q(q-1),\qquad|S_1|=q-2,\qquad|S_2|=|S_3|=q+1.
\]
The four sets are pairwise disjoint. Indeed, $S_0$ has nonzero first coordinate, whereas $S_1$ has first coordinate zero; moreover the third coordinates distinguish $S_2$ and $S_3$ from the preceding sets, taking the values $-1$ and $0$, respectively, while the values $0,-1$ are excluded in the definition of $S_1$. Hence,
\[
|S|=q(q-1)+(q-2)+2(q+1)=q^2+2q,
\]
and therefore
\[
|\Omega|=|S|+1=(q+1)^2.
\]

Since $\Omega$ is a partial ovoid, no two distinct points $P_x,P_y\in\Omega\setminus\{P_\infty\}$ are orthogonal. By Lemma~\ref{lem:affine-chart}, this is equivalent to
\[
Q(x-y)\neq0
\]
for all distinct $x,y\in S$.
\end{proof}

The affine description above can now be inserted directly into the quotient construction. Therefore, we can now prove Theorem \ref{thm:explicit-CP-code}.

\begin{proof}[Proof of Theorem \ref{thm:explicit-CP-code}]
By Proposition~\ref{prop:CP-affine-form}, $|S|=q^2+2q$ and
\[
Q(x-y)\neq0 \qquad\text{for all distinct }x,y\in S.
\]
Therefore, Proposition~\ref{prop:quotient-construction} shows that $\pi^{-1}(S)$ is a symmetric $2$-code.

By the definition of the quotient map,
\[
\pi\left(\begin{pmatrix}
a&b&Z\\
b&c&T\\
Z&T&U
\end{pmatrix}\right)=(c-\delta a,Z,T,U),
\]
and hence
\[
\pi^{-1}(R,Z,T,U)=\left\{\begin{pmatrix}
a&b&Z\\
b&\delta a+R&T\\
Z&T&U
\end{pmatrix}\,:\,a,b\in\F_q\right\}.
\]
Consequently, $\pi^{-1}(S)=\mathcal C_{\delta,\lambda}$.

Since the kernel of $\pi$ has size $q^2$, we obtain
\[
|\mathcal C_{\delta,\lambda}|=q^2|S|=q^2(q^2+2q)=q^4+2q^3.
\]

The minimum rank distance is exactly $2$.  Indeed, the quotient construction gives minimum distance at least $2$, while two distinct matrices in the same fibre differ by a nonzero element of $H_\delta$, and every nonzero element of $H_\delta$ has rank $2$.

Finally,
\[
|\mathcal C_{\delta,\lambda}|=q^4+2q^3>q^4=M_{\mathrm{add}}^{\mathrm{sym}}(3,2),
\]
so $\mathcal C_{\delta,\lambda}$ is non-additive.
\end{proof}

\section{A recursive construction in higher dimension}

The constructions above produce large symmetric $2$-codes in $\Sym_3(q)$. We now show that symmetric $2$-codes in $\Sym_{n_1}(q)$ and $\Sym_{n_2}(q)$ can be combined to produce a symmetric $2$-code in $\Sym_{n_1+n_2}(q)$.  The construction is elementary and uses the off-diagonal block freely.

\begin{theorem}
\label{thm:block-product}
Let $\mathcal C\subseteq\Sym_{n_1}(q),\mathcal D\subseteq\Sym_{n_2}(q)$
be symmetric $2$-codes. Then
\begin{equation}
\mathcal C\star\mathcal D = \left\{\begin{pmatrix}
A&X\\
X^T&B
\end{pmatrix}\,:\,A\in\mathcal C,\;B\in\mathcal D,\;X\in\F_q^{\,n_1\times n_2}\right\}\subseteq\Sym_{n_1+n_2}(q)
\label{eq:block-product}
\end{equation}
is a symmetric $2$-code.  Moreover,
\begin{equation*}
|\mathcal C\star\mathcal D| = q^{n_1n_2}|\mathcal C||\mathcal D|.
\end{equation*}
Consequently,
\begin{equation*}
A_q^{\rm sym}(n_1+n_2,2)
\ge
q^{n_1n_2}
A_q^{\rm sym}(n_1,2)
A_q^{\rm sym}(n_2,2).
\end{equation*}
\end{theorem}

\begin{proof}
Take two distinct matrices in $\mathcal C\star\mathcal D$ and write their difference as
\[
M=\begin{pmatrix}
A-A'&X-X'\\
(X-X')^T&B-B'
\end{pmatrix}.
\]
We prove $\rk M\ge2$.

Suppose first that $A\neq A'$ or$B\neq B'$. Either $A-A'$ or $B-B'$ has rank at least two, because $\mathcal C$ and $\mathcal D$ are $2$-codes.  Since the rank of a matrix is at least the rank of any of its submatrices, $\rk M\ge2$.

Instead, if $A=A'$ and $B=B'$. Since the two matrices are different we must have $X\neq X'$. Hence, put $Y=X-X'\neq0$, then there exists an entry $Y_{ij}\neq0$.  Consider the principal $2\times2$ submatrix of $M$ obtained from the $i$-th coordinate of the first block and the $j$-th coordinate of the second block. It has the form
\[
\begin{pmatrix}
0&Y_{ij}\\
Y_{ij}&0
\end{pmatrix}
\]
and its rank is $2$, therefore $\rk M\ge2$. Thus, every nonzero difference of two elements of $\mathcal C\star\mathcal D$ has rank at least two.

The representation in \eqref{eq:block-product} is unique, and there are $q^{n_1n_2}$ possible matrices $X$. Therefore
\[
|\mathcal C\star\mathcal D|=q^{n_1n_2}|\mathcal C||\mathcal D|,
\]
which proves the result.
\end{proof}

\begin{remark}
The presence of the arbitrary off-diagonal block $X$ is essential for the strength of Theorem~\ref{thm:block-product}.  A block-diagonal product would only have size $|\mathcal C||\mathcal D|$, whereas \eqref{eq:block-product} gains the additional factor $q^{n_1n_2}$. The reason why this is possible is that every nonzero change in $X$ already creates a nonsingular $2\times2$ minor.
\end{remark}

The theorem immediately gives recursive constructions by choosing a small symmetric $2$-code for one of the two factors.

\begin{proposition}
\label{prop:H}
There exist $\alpha,\beta\in\F_q$ such that

\begin{equation*}
\mathcal H_{\alpha,\beta}=\left\{
\begin{pmatrix}
a&b\\
b&\alpha a+\beta b
\end{pmatrix}\,:\,a,b\in\F_q\right\}\subseteq\Sym_2(q)
\end{equation*}

is an additive symmetric $2$-code of cardinality $q^2$.  
\end{proposition}

\begin{proof}
Choose $\alpha,\beta\in\F_q$ such that $f(x)=x^2-\beta x-\alpha$ is irreducible over $\F_q$; such polynomials exist for every $q$. For $(a,b)\neq(0,0)$,
\[
\det\begin{pmatrix} 
a&b\\ 
b&\alpha a+\beta b
\end{pmatrix}=-(b^2-\beta ab-\alpha a^2).
\]
If $a=0$ this equals $-b^2\neq0$; if $a\neq0$ it equals $-a^2 f(b/a)\neq0$. Hence, every nonzero element of $\mathcal H_{\alpha,\beta}$ is nonsingular.
\end{proof}

\begin{remark}
The existence of an additive symmetric $2$-code in $\Sym_2(q)$ of cardinality $q^2$ also follows directly from the field extension $\F_{q^2}/\F_q$. Indeed, fixing an $\F_q$-basis $\{e_1,e_2\}$ of $\F_{q^2}$, consider the set
\[
\mathcal H =\left\{\left(\Tr_{\F_{q^2}/\F_q}(x e_i e_j)\right)_{i,j=1}^2\,:\,\ x\in\F_{q^2}\right\}.
\]
Since the trace pairing is nondegenerate, every nonzero matrix in $\mathcal H$ is nonsingular. Thus, $\mathcal H$ is an additive symmetric $2$-code of cardinality $q^2$. Proposition~\ref{prop:H} provides a particularly simple explicit representation of such a code.
\end{remark}

Applying Theorem~\ref{thm:block-product} with $\mathcal D=\mathcal H_{\alpha,\beta}$ yields the following recursion.

\begin{corollary}
\label{cor:two-dimensional-recursion}

For every $n\ge1$,
\begin{equation*}
A_q^{\rm sym}(n+2,2)\ge q^{2n+2}A_q^{\rm sym}(n,2).
\end{equation*}
More precisely, if $\mathcal C_n\subseteq\Sym_n(q)$ is a symmetric $2$-code, then
\begin{equation*}
\mathcal C_{n+2}=\left\{
\begin{pmatrix}
A&X\\
X^T&H
\end{pmatrix}\,:\,A\in\mathcal C_n,\;X\in\F_q^{\,n\times2},\;H\in\mathcal H_{\alpha,\beta}\right\}
\end{equation*}
is a symmetric $2$-code with $|\mathcal C_{n+2}|=q^{2n+2}|\mathcal C_n|$.
\end{corollary}

\begin{proof}
Apply Theorem~\ref{thm:block-product} with $n_2=2$ and $|\mathcal H_{\alpha,\beta}|=q^2$, then
\[
|\mathcal C_{n+2}|=q^{2n}|\mathcal C_n|q^2=q^{2n+2}|\mathcal C_n|.
\]
\end{proof}

We can now propagate any improvement obtained in dimension three to all odd dimensions.

\begin{theorem}
\label{thm:propagation-seed}
Suppose that there exists a symmetric $2$-code $\mathcal C_3\subseteq\Sym_3(q)$ of size $M$.  Then, for every odd integer $n\ge3$,
\begin{equation*}
A_q^{\rm sym}(n,2)\ge q^{(n^2-9)/2}M.
\end{equation*}
\end{theorem}

\begin{proof}
Write $n=3+2r$ and, starting with $\mathcal C_3$, apply Corollary~\ref{cor:two-dimensional-recursion} successively in dimensions $3,\ 5,\ 7,\ldots,\ n-2$. The total power of $q$ introduced by the recursion is
\[
\sum_{j=0}^{r-1}
\bigl(2(3+2j)+2\bigr)
=
\sum_{j=0}^{r-1}(8+4j)=2r^2+6r=\frac{n^2-9}{2},
\]
since $r=\frac{n-3}{2}$. Therefore, the recursively constructed code has size $q^{(n^2-9)/2}M$.
\end{proof}

If we apply this iteration to the Cossidente-Marino-Pavese code we obtain the following result.
\begin{corollary}
\label{cor:all-symmetric}
    Let $q,n\geq 3$ with $n$ odd. Then
    \begin{equation*}
        A_q^{\rm sym}(n,2) \ge q^{(n^2-1)/2}+q^{(n^2-3)/2}+q^{(n^2-9)/2}.
    \end{equation*}
\end{corollary}
\begin{proof}
Since $q\geq3$ there exists a Cossidente-Marino-Pavese code $\mathcal{C}_{\rm CMP}\subseteq \Sym_3(q)$. Such a code is a symmetric $2$-code with cardinality $|\mathcal{C}_{\rm CMP}|=q^4+q^3+1$, therefore, applying Theorem \ref{thm:propagation-seed}, we get
    \begin{equation*}
        A_q^{\rm sym}(n,2) \ge q^{(n^2-9)/2}(q^4+q^3+1)= q^{(n^2-1)/2}+q^{(n^2-3)/2}+q^{(n^2-9)/2}.
    \end{equation*}
\end{proof}

We now apply Theorem \ref{thm:propagation-seed} to the codes obtained from partial ovoids.

\begin{corollary}
\label{cor:partial-ovoid-propagation}
Let $q$ be odd and suppose that $Q^-(5,q)$ contains a partial ovoid of size $m$.  Then, for every odd $n\ge3$,
\begin{equation}
A_q^{\rm sym}(n,2) \ge q^{(n^2-5)/2}(m-1).
\label{eq:partial-ovoid-all-odd}
\end{equation}
\end{corollary}

\begin{proof}
By Theorem~\ref{prop:partial-ovoid-input}, the partial ovoid gives a symmetric $2$-code in $\Sym_3(q)$ of size $M=q^2(m-1)$. Theorem~\ref{thm:propagation-seed} therefore gives
\[
A_q^{\rm sym}(n,2)\ge q^{(n^2-9)/2}q^2(m-1)=q^{(n^2-5)/2}(m-1).
\]
\end{proof}

Thus, every improvement in the size of a partial ovoid of $Q^-(5,q)$ automatically produces an improvement in every odd matrix dimension.

Using the Cossidente--Pavese partial ovoid of size $m=(q+1)^2$ gives a particularly simple uniform family.

\begin{corollary}
\label{cor:all-odd-symmetric}
Let $q$ and $n\ge3$ be odd.  Then
\begin{equation*}
A_q^{\rm sym}(n,2)\ge q^{(n^2-1)/2}+2q^{(n^2-3)/2}.
\end{equation*}
\end{corollary}

\begin{proof}
Since $m=(q+1)^2$, we have $m-1=q^2+2q=q(q+2)$.
Substituting this into
\eqref{eq:partial-ovoid-all-odd} gives
\[
A_q^{\rm sym}(n,2)\ge q^{(n^2-5)/2}q(q+2) = q^{(n^2-3)/2}(q+2)=q^{(n^2-1)/2}+2q^{(n^2-3)/2}.
\]
\end{proof}

The largest cardinality of an additive symmetric $2$-code in
$\Sym_n(q)$ is
$q^{\lfloor n^2/2\rfloor}$;
see, for instance, \cite{Schmidt2015}.  In particular, when $n$ is
odd, the additive bound is
\begin{equation*}
M_{\rm add}^{\rm sym}(n,2)=q^{(n^2-1)/2}.
\end{equation*}

Corollary \ref{cor:all-symmetric} and Corollary~\ref{cor:all-odd-symmetric} therefore exceed the additive bound in every odd dimension $n\ge3$.  More precisely, doing the ratio between the bound obtained, respectively, in Corollary \ref{cor:all-symmetric} and Corollary~\ref{cor:all-odd-symmetric} and the known value for $M_{\rm add}^{\rm sym}(n,2)$, we get
\begin{equation*}
\frac{q^{(n^2-9)/2}(q^4+q^3+1)}{q^{(n^2-1)/2}}=1+\frac{1}{q}+\frac{1}{q^4}, \qquad q\geq 3
\end{equation*}
and
\begin{equation*}
\frac{q^{(n^2-3)/2}(q+2)}{q^{(n^2-1)/2}}=1+\frac{2}{q}, \qquad q \text{ odd}.
\end{equation*}
Thus, the multiplicative improvement obtained in dimension
three is preserved by the recursion.

\section{Concluding remarks and open problems}
\label{sec:concluding}

In this paper, we have investigated non-additive symmetric rank-metric codes with minimum distance two.  Our main construction starts with a two-dimensional constant-rank subspace $H_\delta\leq\Sym_3(q)$ and passes to the quotient $\Sym_3(q)/H_\delta$. In this quotient, the obstruction to minimum distance two is described by the isotropic cone of a nondegenerate quadratic form on $\F_q^4$. By embedding this quadratic space in an elliptic quadric $Q^-(5,q)$, sets with
anisotropic pairwise differences are naturally obtained from partial ovoids.

This correspondence gives, for every partial ovoid of $Q^-(5,q)$ of size $m$, a symmetric $2$-code in $\Sym_3(q)$ of cardinality $q^2(m-1)$.

In particular, the Cossidente--Pavese construction yields, for every odd prime power $q$, $A_q^{\rm sym}(3,2)\ge q^4+2q^3$, and hence produces symmetric rank-metric codes which are strictly larger than the largest additive codes with the same parameters. We have also given an explicit matrix description of this family.

Moreover, we propagated the improvements obtained in small dimension to larger dimensions using the block construction of Theorem~\ref{thm:block-product}. For example, combining it with the partial-ovoid construction gives, for every odd $n,q\geq3$,

$$
A_q^{\rm sym}(n,2)\ge q^{(n^2-1)/2}+2q^{(n^2-3)/2},
$$

again exceeding the additive bound

$$
M_{\rm add}^{\rm sym}(n,2)=q^{(n^2-1)/2}.
$$

Several natural questions remain open.  The most immediate one concerns the three-dimensional extremal problem.

\begin{problem}
Determine the exact value of
$A_q^{\rm sym}(3,2)$.
\end{problem}

For odd $q$, our construction and the upper bound \eqref{eq:CMP-sym-upper} give
\[
q^4+2q^3\leq A_q^{\rm sym}(3,2) \leq 1+\frac{q(q^2-1)(q^2+q+1)}{2}.
\]
Thus, a substantial gap remains between the best known lower and upper bounds.

The connection with partial ovoids provides a natural geometric approach to this problem.  If $\mu^-(5,q)$ denotes the maximum size of a partial ovoid of $Q^-(5,q)$, then our construction gives
\[
A_q^{\rm sym}(3,2)\geq q^2\bigl(\mu^-(5,q)-1\bigr).
\]
It would therefore be interesting to understand whether improvements can be obtained through larger partial ovoids, or whether extremal symmetric $2$-codes require constructions that do not arise from the quotient model developed here.

A second question concerns higher dimensions.  Our block construction propagates the three-dimensional improvement to every odd dimension, but it does not exploit the geometry of $\Sym_n(q)$ directly.

\begin{problem}
For odd $n$ and $q\geq 3$ construct symmetric $2$-codes in $\Sym_n(q)$ that improve upon
\[
A_q^{\rm sym}(n,2)\geq q^{(n^2-1)/2}+2q^{(n^2-3)/2},\qquad \text{when $q$ is odd},
\]
\begin{equation*}
    A_q^{\rm sym}(n,2) \ge q^{(n^2-1)/2}+q^{(n^2-3)/2}+q^{(n^2-9)/2},\qquad \text{when $q$ is even.}
\end{equation*}
\end{problem}

\section*{Declaration of Generative AI and AI-assisted Technologies in the Writing Process}

During the preparation of this work, the authors used GPT-5.6 Sol (OpenAI) to produce a first draft of the proofs in Sections~3--6 from the authors' notes, and to assist with the computation of the explicit coordinates of the Cossidente--Pavese partial ovoid in Section~5.
After using this tool, the authors reviewed and edited the content, independently verified all mathematical arguments, computations and bibliographic references, and take full responsibility for the content of the published article.

\section*{Acknowledgments}
Marta Messia is very grateful for the hospitality of the Department of Mathematics and Physics, Università della Campania Luigi Vanvitelli, Caserta, where she spent two weeks as a visiting researcher.

The research of Marta Messia is supported by FWO research project ``Finite Geometry Applications to Coding Theory and Cryptography’’ (FWO grant nr. G0AJI25N). This research was also carried out as part of the Scientific Research Network ``Graphs, Association schemes and Geometries: structures, algorithms and computation” of the Research Foundation Flanders (FWO), (FWO network nr. W003324N).

The research of Ferdinando Zullo was partially supported by the Italian National Group for Algebraic and Geometric Structures and their Applications (GNSAGA - INdAM).
This work was partially written during the Opera 2026 conference in Bordeaux, and so we acknowledge the support from the International Research Laboratory LYSM in partnership between CNRS and INdAM.

\end{document}